\documentclass[a4paper,fleqn,preprintnumbers]{elsarticle}
\usepackage{booktabs}
\usepackage{enumitem}
\usepackage{tikz}
\usepackage{graphicx}
\usepackage{soul}
\usepackage{nameref}
\usepackage{bm}
\usepackage{amssymb} 
\usepackage[english]{babel}
\usepackage{hyperref}
\usepackage{amsmath}
\newtheorem{thm}{Theorem}
\newtheorem{lem}{Lemma}

\def\tsc#1{\csdef{#1}{\textsc{\lowercase{#1}}\xspace}}
\tsc{WGM}
\tsc{QE}

\begin{document}
\let\WriteBookmarks\relax
\def\floatpagepagefraction{1}
\def\textpagefraction{.001}



\title{Coprime networks of the composite numbers: \\ pseudo-randomness and synchronizability}  

\author[1]{Md Rahil Miraj}


\fnmark[1]

\ead{rahilmiraj@gmail.com}


\affiliation[1]{organization={Theoretical Statistics and Mathematics Unit, Indian Statistical Institute},
            addressline={8th Mile, Mysore Rd, RVCE Post}, 
            city={Bengaluru},
            postcode={560059}, 
            state={Karnataka},
            country={India}}
            
\author[2]{Dibakar Ghosh}

\fnmark[2]

\ead{dibakar@isical.ac.in}


\ead[url]{}


\affiliation[2]{organization={Physics and Applied Mathematics Unit, Indian Statistical Institute},
            addressline={203 B.T. Road}, 
            city={Kolkata},
            postcode={700108}, 
            state={West Bengal},
            country={India}}

\author[3]{Chittaranjan Hens}

\fnmark[3]

\ead{chittaranjanhens@gmail.com}

\affiliation[3]{organization={Center for Computational Natural Science and Bioinformatics, International Institute of Informational Technology},
            addressline={Gachibowli}, 
            city={Hyderabad},
            postcode={500032}, 
            state={Telengana},
            country={India}}
        
\maketitle

In this paper, we propose a network whose nodes are labeled by the composite numbers and two nodes are connected by an undirected link if they are relatively prime to each other. As the size of the network increases, the network will be connected whenever the largest possible node index $n\geq 49$. To investigate how the nodes are connected, we analytically describe that the link density saturates to $6/\pi^2$, whereas the average degree increases linearly with slope $6/\pi^2$ with the size of the network. To investigate how the neighbors of the nodes are connected to each other, we find the shortest path length will be at most 3 for $49\leq n\leq 288$ and it is at most 2 for $n\geq 289$. We also derive an analytic expression for the local clustering coefficients of the nodes, which quantifies how close the neighbors of a node to form a triangle. We also provide an expression for the number of $r$-length labeled cycles, which indicates the existence of a cycle of length at most $O(\log n)$. Finally, we show that this graph sequence is actually a sequence of weakly pseudo-random graphs. We numerically verify our observed analytical results.  As a possible application, we have observed less synchronizability (the ratio of the largest and smallest positive eigenvalue of the Laplacian matrix is high) as compared to Erd\H{o}s–R\'{e}nyi random network and Barab\'{a}si-Albert network. This unusual observation is  consistent with the prolonged transient behaviors of  ecological and predator-prey networks which can easily avoid  the global synchronization.\\



\section{Introduction}
The study of graphs with intricate connections between the items is known as a complex network \cite{Newman_Network_Book,albert2002statistical,pastor2003statistical,boccaletti2006complex}. The field of network science examines a wider range of systems found in the real world \cite{arenas2008synchronization,barrat2008dynamical,Cohen_complex_book,battiston2021physics,hens2019spatiotemporal,gao2016universal,barzel2013universality,meena2023emergent},\cite{artime2024robustness}. For instance, the study of network science can be used to map out emerging behavior in social networks \cite{borgatti2009network}, the proper functionality of the technological systems (such as the power grid) \cite{menck2013basin}, the spreading mechanism of diseases among communities \cite{wang2019coevolution,pastor2015epidemic}, the stability of prey-predator systems \cite{donohue2016navigating}, neuronal activities \cite{bullmore2012economy}, signal propagation\cite{ji2023signal} and the interaction mechanism in subcellular systems \cite{barabasi2004network}.  

The study of different mathematical properties of some sets of numbers is now frequently done using network science \cite{shekatkar2015divisibility,garcia2014complex}.  Nodes from those sets are used to build the networks, which connect each pair of them by directed or undirected links whenever they fulfill a particular relation.
Many such interesting networks have been found in the literature.
For instance, divisibility networks of natural numbers following usual increasing sequential order \cite{shekatkar2015divisibility}  or by taking nodes from the Pascal matrix and constructing a divisibility network with those nodes \cite{solares2020divisibility} are recently developed. Also in 
 \cite{chandra2005small}, assuming that Goldbach conjecture is true, the authors have broken the even numbers by the sum of two primes and constructed a network of those prime numbers by connecting links between two primes depending on the corresponding even numbers.  Sometimes, it is also possible to find analytic expressions of the measures like average degree, link density, clustering coefficients etc. and by constructing networks of different size one may check how fast it agrees with the theoretical results. In 
   \cite{shekatkar2015divisibility} the authors have shown numerically that the divisibility network of natural numbers follows a scale-free degree distribution. They have shown analytically that the global clustering coefficient of this network decays to zero whereas the average degree increases logarithmically as the size of the network increases.

A random graph of $n$ vertices is a graph, whose edges between two vertices are drawn randomly. On the other hand, a pseudo random graph is a graph that behaves like a random graph of same number of vertices and same edge density with a high probability. In order to give a quantitative measure, Andrew Thomason first introduced the notion of $(p,\alpha)$-jumbled graph in 1987 \cite{thomason1987pseudo}. He defined a graph $G=(V,E)$ is said to be $(p,\alpha)$-jumbled graph for  real numbers $p$ and $\alpha$ with $0<p<1\leq \alpha$, if for every induced subgraph $H$ of $G$, we get 
$$\left|E(H)-p\binom{|H|}{2} \right|\leq \alpha|H|.$$
A $(p, O(\sqrt{np} ))$-jumbled graph is one of the best possible pseudo-random graph, as  Erd\H{o}s–R\'{e}nyi random graph G(n,p) is almost surely $(p, O(\sqrt{np} ))$-jumbled \cite{krivelevich2006pseudo}. 

Since, the largest eigenvalue of the adjacency matrix of a $d$-regular graph is $d$, Alon 
  \cite{krivelevich2006pseudo,alon1985lambda1,alon1986eigenvalues} defined by $(n, d, \lambda)$-graph, a $d$-regular graph of $n$ vertices and eigenvalues $d=\lambda_1>\lambda_2\geq...\geq\lambda_n$ such that $|\lambda_i|\leq \lambda$ for all $i=2,3,...,n$. It has been shown that $(n, d, \lambda)$-graph with $\lambda$ much smaller than $d$ have certain pseudo-random properties 
 \cite{alon1985lambda1,alon1986eigenvalues,alon2016probabilistic}. And it has also been found in literature that a $(n, d, \lambda)$-graph is $(d/n,\lambda)$-jumbled \cite{krivelevich2006pseudo}. Apart from the pseudo-random graphs of the type $(p,\alpha)$-jumbled graphs, in general, in a pseudo-random graph $G$ with any red/blue-coloring of edges such that the proportion of red colored edges is $r$, it has been found that there exists a Hamiltonian cycle having the proportion of red edges is close to $r$ \cite{kuhn2006multicolored}. Pancyclicity, i.e, the existence of cycles of all possible lengths in pseudo-random graphs has also been explored in literature \cite{krivelevich2010resilient}.

By slightly relaxing the conditions of $(p,\alpha)$-jumbled graph, we define a sequence of graphs $\{G_n\}$ to be weakly pseudo-random \cite{krivelevich2006pseudo} if for all subsets $U\subseteq V(G_n)$, we have
$$\left|E(U)-p\binom{|U|}{2}\right|=o(pn^2),$$
where  $p = p(n)$ is a parameter, which is typically the edge density of graphs in the sequence. Chung, Graham and Wilson first defined such a graph with $p=1/2$ and named quasi-random graphs \cite{chung1989quasi} in 1989 and also provided several equivalent statements, which are actually same as above. 

In this study, we construct a network with composite numbers as nodes and if two of them are relatively prime to one another, we connect them with an undirected link. In this case, we are looking for answers to certain pertinent concerns about the structural characteristics of graphs from the standpoint of graph theory. For instance, how the nodes are linked together and how many nearby nodes may a node have? 
We also look for the diameter of the graph and try to understand how the average shortest path length or the diameter of the graph depend on the size of the network. We also want to know how a node's neighbors are connected between themselves. This property is measured by the local clustering coefficient of a node, which is defined as the proportion of the number of neighbors of a node being connected to each other to $\binom{d}{2}$, where $d$ is the degree of that node \cite{barabasi2013network}(sec-2.10).

On the other hand, we revisit the old question of composite number theory: what is the probability that two randomly chosen composite numbers to be coprime to each other, which will be answered in the subsequent sections by analytically calculating the average degree to be $\frac{6n}{{\pi}^2}$ and link density $\frac{6}{{\pi}^2}$.
Also the existence of labeled cycles of a given length is another question that one may ask. Here we not only have shown the existence of such a labeled cycle of some given length but also calculate an expression for it, which gives us a glimpse of the existence of a cycle of length at most $O(\log n)$ for large $n$. Therefore the emergence of such dense graph motivates us to identify or categorize the structural pattern of the graph, which type of graph it is. So, to answer this question, we are able to show that the graph sequence namely coprime networks of composite numbers are actually weakly pseudo-random with $p=\frac{6}{{\pi}^2}$. 
As an application we have studied the eigenvalues of the Laplacian matrix of the constructed network and infer that the synchronizability  in the coprime network of composite numbers is quite  less as compared to Erd\H{o}s–R\'{e}nyi random network and Barab\'{a}si-Albert network \cite{arenas2008synchronization,barahona2002synchronization,pecora1998master,nishikawa2003heterogeneity,osipov2007synchronization}. This less synchronizability is noticed in ecological and predator-prey networks. Synchronization is a natural phenomenon that is mostly observed in a population of dynamically interacting elements. Synchronization processes play a very
important role in many different contexts such as biology, ecology, climatology, sociology, technology \cite{pikovsky2001synchronization, osipov2007synchronization}. It can be seen when oscillating elements are constrained to interact in a complex network topology\cite{arenas2008synchronization}. Note that, two aspects in particular affect the network synchronizability entirely: the first is the synchronized region associated with the node dynamics and the second has to do with the eigenvalues of the network structural matrix. In \cite{dai2020discontinuous}, the authors demonstrated that explosive synchronization occurs at all even dimensions when Kuramoto oscillators are used in higher dimensions. In this paper, we have significantly used the eigenvalues of the Laplacian matrix of our constructed network to infer about the synchronizability.
\begin{figure}[htbp]
    \centering
        \includegraphics[scale=0.4]{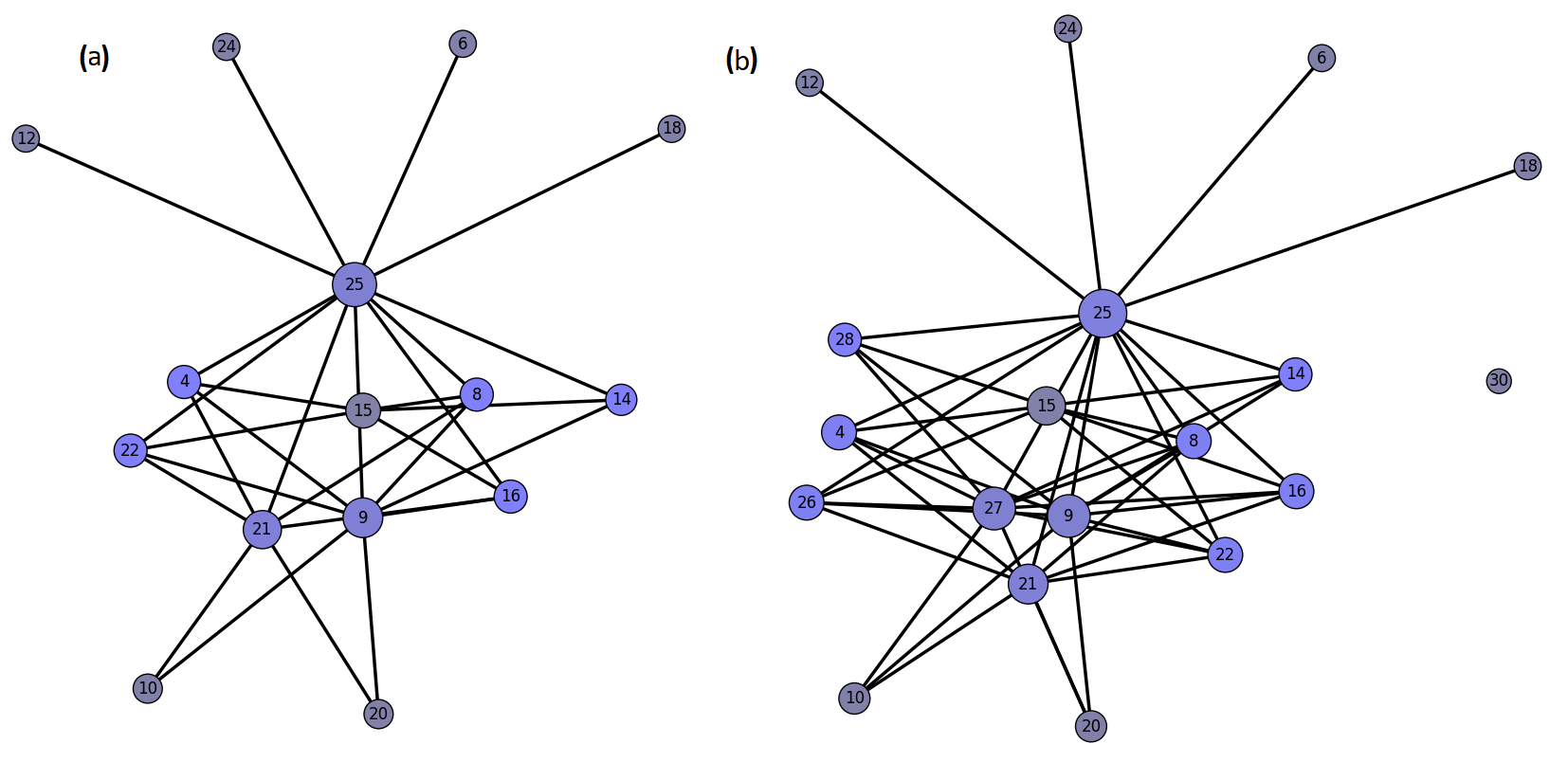}
       \caption{\textbf{Coprime networks of composite numbers:} (a) $n=25$, (b) $n=30$. Here, the size of each node is proportional to its degree and the color of each node is given according to its clustering coefficient with lighter blue nodes as nodes with higher value of local clustering. For $n=30$, the network is disconnected. Later we will show that the constructed network will be connected if $n\geq 49$.}
\label{Fig:1}
\end{figure}

\section{Construction of the Network}
To construct the proposed network, we chose the composite numbers less than or equal to $n$ as nodes and connect two of them by an undirected edge if they are coprime to each other, i.e., for $u,v\in \mathbf{V},~u\sim v$, if and only if $gcd(u,v)=1$, where $\mathbf{V}$ is the set of nodes, i.e., the set of composite numbers less than or equal to $n$. Since, we have a natural ordering in the composite numbers, so at each step we include the next composite number as a new node which is greater than the largest node in the existing network and connect it to those nodes which are coprime to it. Thus we construct the network in a growing manner. In the constructed network, we have an identity of each node as a composite number and we can study different properties of the network as a function of the node index. This is not possible with other networks, where there is no specific node index.

\section{Properties of the Network}
Now, we see the schematic representation of the networks for $n=25$ and $n=30$. In this section, we prove a few important properties of the constructed network. First we are eager to look analytically what is the least composite number after which we get a connected network always. In the next theorem, we will show that the occurrence of isolated node is not random, rather as we increase $n$, after a certain point we will never have an isolated node.

\begin{thm} \label{Theorem 1}
There doesn't exist any isolated node in the network constructed as above with largest possible node index $n\geq 49$.
\end{thm}

\textit{Proof.}
We can see that if there is no composite number up to $n$ that is coprime to some composite number $k$, then the node $k$ will be isolated in that graph. If $k={p_1}^{a_1}{p_2}^{a_2}\cdots {p_r}^{a_r}$, where $p_i$'s are distinct primes such that $p_1=2$, $p_2=3$ and so on and $a_i$'s are non-negative integers, except $a_r\neq 0$. Now if some $a_i$ is $0$, then $k$ will be connected to ${p_i}^2$, so it will be enough to take $a_i=1$ for all $i=1,2,\cdots,r$. Then the least node greater than $k$ that will be connected to $k$ is the square of the least prime number $p\neq p_i$ for $i=1,2,\cdots,r$, i.e., for $k=p_1p_2\cdots p_r$, the smallest node that will be connected to $k$ is $p_{r+1}^2$, but if $p_{r+1}^2>n$, then $k$ will be isolated. 

So, for some $n$, if $p_{r+1}^2<p_1p_2\cdots p_r= k_n\leq n$, but $p_1p_2\cdots p_{r+1}>n$, then if $p_1p_2.\cdots p_rp_{r+1}>n+1$, we will get $k_{n+1}=k_n=p_1p_2\cdots p_r$, which will be connected to $p_{r+1}^2$ and if $p_1p_2\cdots p_rp_{r+1}=n+1$, then by Lemma \ref{Lemma 6} in the appendix with $t=2$ and $s=r$, we get $n+1=k_{n+1}=p_1p_2\cdots p_rp_{r+1}>p_{r+2}^2$. So, $k_{n+1}$ will be connected to $p_{r+2}^2$ for the latter case.

But the smallest $r$ satisfying $p_1p_2\cdots p_r> p_{r+1}^2$ is $r=4$. So, by induction, we can conclude that we will always have $p_1p_2\cdots p_r> p_{r+1}^2$ for $r\geq 4$. Thus, for $n\geq p_4^2=7^2=49$, there will be no isolated point.

\hfill \ensuremath{\Box}

Now, we can see that in a particular network with largest possible node index $n$, the total number of nodes is the total number of composite numbers less than or equal to $n$.
\begin{thm}\label{Theorem 2}
The number of nodes in the proposed network with largest possible node index $n$ is $N(n)=(n-\pi(n)-1)$, where $\pi(n)$ is the total number of primes less than or equal to $n$.
\end{thm}

Now, we will calculate the total number of edges of a network with greatest possible node index $n$.

\begin{thm} \label{Theorem 3}
The number of edges in a particular network with largest possible node index $n$ is
\[E(n)=\frac{3}{{\pi}^2}n^2+n \log \log n-n\pi(n)+\frac{\pi(n)(\pi(n)+1)}{2}+O(n \log n).\]
\end{thm}

\textit{Proof.}
 If we add nodes by increasing order starting from 4 and connect them by edges accordingly at each step, the number of new added edges after adding the $k^{th}$ node is
\[E(k)-E(k-1)=
\begin{cases}
  \phi(k)-\pi(k)+w(k)-1, & \text{if $k$ is $composite$} \\
  0, & \text{otherwise},
\end{cases}\]
where $\phi(k)$ is the total number of positive integers less than $k$ which are coprime to $k$ 
 \cite{apostol1998introduction} and $w(k)$ is the total number of distinct prime factors of $k$.

Therefore, the telescoping sum gives
\[E(n)-E(1)=\sum_{\substack{k=2 \\ k ~\text{composite}}}^{n}\Big[ \phi(k)-\pi(k)+w(k)-1 \Big] \]
gives,
\begin{equation}\label{1}
E(n)=\sum_{\substack{k=2 \\ k \neq \text{p}}}^{n}\Big[ \phi(k)-\pi(k)+w(k)-1 \Big],
\end{equation}
as $E(1)=0$, where $p$ is prime number.

We can write it as,
\begin{align*}
E(n) & =\sum_{k=2}^{n}\Big[\phi(k)-\pi(k)+w(k)-1\Big]-\sum_{p\leq n}\Big[\phi(p)-\pi(p)+w(p)-1\Big]\\
& =\sum_{k=2}^{n}\Big[\phi(k)-\pi(k)+w(k)-1\Big]-\sum_{i=1}^{m}\Big[\phi(p_i)-\pi(p_i)\Big],   
\end{align*}

where $m=\pi(n)$ and $w(p_i)=1$ for all $i=1,2,...,m$, where $p_i$ is the i-th prime.
\begin{align*}
& =\sum_{k=2}^{n}\Big[\phi(k)-\pi(k)+w(k)-1\Big]-\sum_{i=1}^{m}\Big[(p_i-1)-(i-1)\Big]\\
& =\sum_{k=2}^{n}\Big[\phi(k)-\pi(k)+w(k)-1\Big]-\sum_{i=1}^{m}(p_i-i).
\end{align*}

Therefore,
\begin{equation}\label{2}
E(n)=\sum_{k=2}^{n}\phi(k)-\sum_{k=2}^{n}\pi(k)+\sum_{k=2}^{n}w(k)-(n-1)-\sum_{i=1}^{m}p_i+\frac{m(m+1)}{2}.
\end{equation}
Now, using Lemmas \ref{Lemma 1}, \ref{Lemma 2}, \ref{Lemma 5} in the appendix and from (\ref{2}), we get
$$E(n)=\Big[\frac{3}{{\pi}^2}n^2+O(n \log n)\Big]-\Big[mn -\sum_{i=1}^{m}p_i\Big]+\Big[n \log \log n+B_1n+O(n)\Big]-(n-1)-\sum_{i=1}^{m}p_i+\frac{m(m+1)}{2}$$
$$=\frac{3}{{\pi}^2}n^2+n \log \log n-nm+\frac{m(m+1)}{2}+O(n \log n).~~~~~~~~~~~~~~~~~~~~~~~~~~$$

\hfill \ensuremath{\Box}

As we add all the composite numbers up to the natural number $n$ in the graph one by one as the nodes and at each step we join any two of them by undirected links if they are relatively prime, we get the desired network for some given $n$. Then, a node $k$ will be connected to all the composite numbers less than or equal to $n$ which are coprime to it and that is the degree of $k$. 

\begin{thm}\label{Theorem 4}
In the proposed network with largest possible node index $n$:\\
i. The degree of a node $k$, $d_k=\phi(n,k)-\pi(n)+w(k)-1$, where $\phi(n,k)$ is the total number of positive integers less than or equal to $n$ which are coprime to $k$.\\
ii. The average degree of the nodes, $\overline{d(n)}=\frac{6}{{\pi}^2}n+O\left(\frac{n}{(\log n)^2}\right),$ i.e., $\overline{d(n)} \sim \frac{6}{{\pi}^2}n$.\\
iii. The maximum degree of a node $k$ is $\Delta\leq n-\sqrt{n}-\frac{n}{\log n}.$\\
iv. The codegree, i.e., the number of common neighbors between two nodes $k$ and $l$ is
$codeg(k,l)=\phi(n,kl)-\pi(n)+w(kl)-1$.
\end{thm}

\textit{Proof.}

i. We have
  $$ d_k =\sum\limits_{\substack{i=2 \\ i\neq \text{prime} \\ (i,k)=1}}^{n}\mathbf{1}=\sum\limits_{\substack{i=2 \\ (i,k)=1}}^{n}\mathbf{1}-\sum\limits_{\substack{p\leq n \\ p\nmid k}}^{n}\mathbf{1}=\phi(n,k)-\pi(n)+w(k)-1.$$

ii. Now,
\begin{align*}
   \overline{d(n)} & =\frac{1}{n-\pi(n)-1}\sum_{k\in \mathbf{V}}d_k=\frac{2E(n)}{n-\pi(n)-1} \\
     & =\frac{2\Big[\frac{3}{{\pi}^2}n^2+n \log \log n-n\pi(n)+\frac{\pi(n)(\pi(n)+1)}{2}+O(n \log n)\Big]}{n-\pi(n)-1} \\
     & \sim \frac{2\Big[\frac{3}{{\pi}^2}n^2+n \log \log n-\frac{n^2}{\log n}+\frac{n}{2\log n}\left(\frac{n}{\log n}+1\right)+O(n \log n) \Big]}{n-\frac{n}{\log n}-1}\\
     & = \frac{2\Big[\frac{3}{{\pi}^2}n+\log \log n-\frac{n}{\log n}+\frac{n}{2(\log n)^2}+\frac{1}{2\log n}+O(\log n) \Big]}{1-\frac{1}{\log n}-\frac{1}{n}} \\
     & =\frac{6}{{\pi}^2}n+O\left(\frac{n}{(\log n)^2}\right).
  \end{align*}
Therefore, $\overline{d(n)} \sim \frac{6}{{\pi}^2}n$.
\begin{figure}[htbp]
    \centering
        \includegraphics[scale=.6]{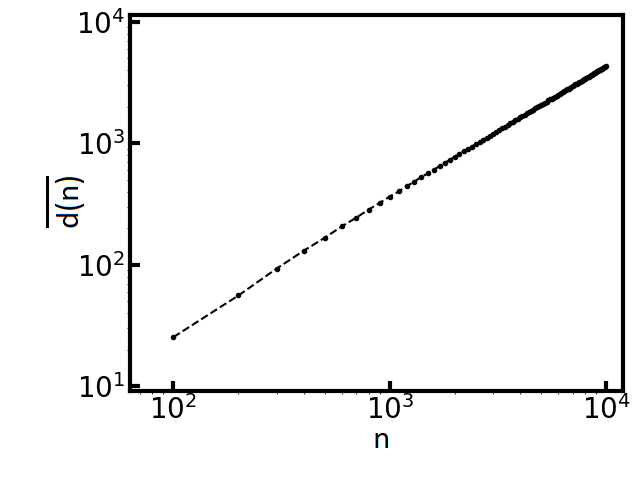}
       \caption{{The average degree $\overline{d(n)}$ is plotted by changing the highest possible node index $n$. Here the average degree increases monotonically by adding the nodes. The data has been plotted using log-log scale.  The scaling exponent is 1.}}
    \label{Fig:2}
    \end{figure}

Theoretical result suggests that for large $n$, the mean degree scales with $n^{\alpha} $, where $\alpha=1$, and the proportionality factor is around $0.60.79(\frac{6}{{\pi}^2})$. The data is plotted in Fig.\ref{Fig:2} with black dots. The data is plotted in log-log scale. Note that, for computational limitation, we have varied network size from 10 to $10^4$. However, the theoretical line fits very well with the numerical data.

Clearly, the network is highly dense, and the prefactor suggests that each node is connected with almost $60.79\%$ nodes of the graph, which further signify that for a given $n$, we can say that the probability of a pair of composite numbers to be coprime is $\sim0.60.79$.

iii. Now, let $r\geq1$ be such that we have $p_r^2\leq n< p_{r+1}^2$. Then we can see that $p_r^2$ will have the maximum degree. So,
$$ \Delta = d_{p_r^2} =\phi(n,p_r^2)-\pi(n)+w(p_r^2)-1 = \phi(n,p_r)-\pi(n) =n-\left\lfloor\frac{n}{p_r}\right\rfloor-\pi(n) \leq n-\lfloor\sqrt{n}\rfloor-\pi(n),$$
since, $p_r^2\leq n$.
Therefore, $$\Delta \sim n-\sqrt{n}-\frac{n}{\log n}.$$

\begin{figure}[htbp]
    \centering
        \includegraphics[scale=.6]{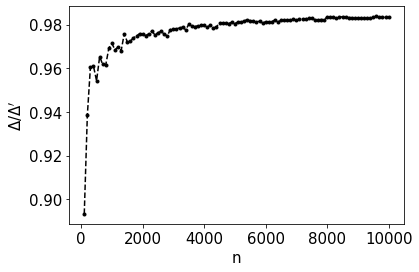}
       \caption{{The ratio of the computed and theoretical max degree($\Delta/\Delta'$) is plotted by changing the highest possible node index $n$. The ratio getting closer to 1 indicates the accuracy of the theoretical maximum of the degree.}}
    \label{Fig:3}
    \end{figure}
    
iv. Again, we have
    
 $$ codeg(k,l) =\sum_{\substack{i=2 \\ i\neq \text{prime} \\ (i,k)=1 \\ (i,l)=1}}^{n}\mathbf{1}= \sum_{\substack{i=2 \\ i\neq \text{prime} \\ (i,kl)=1}}^{n}\mathbf{1} =\sum_{\substack{i=2 \\ (i,kl)=1 }}^{n}\mathbf{1}-\sum_{\substack{p\leq n \\ p\nmid kl}}^{n}\mathbf{1}=\phi(n,kl)-\pi(n)+w(kl)-1.$$

\hfill \ensuremath{\Box}

\begin{thm}\label{Theorem 5}
 The link density in the proposed network with largest possible node index $n$ is
 $$l(n)=\frac{6}{{\pi}^2}+O\left(\frac{\log n}{n}\right),~~\text{i.e.,}~~l(n) \sim \frac{6}{{\pi}^2}.$$
\end{thm}

\textit{Proof.}
  \begin{align*}
    l(n) & =\frac{E(n)}{\binom{N(n)}{2}}=\frac{2E(n)}{(n-\pi(n)-1)(n-\pi(n)-2)} \\
     & =\frac{2\Big[\frac{3}{{\pi}^2}n^2+n \log \log n-n\pi(n)+\frac{\pi(n)(\pi(n)+1)}{2}+O(n \log n)\Big]}{(n-\pi(n)-1)(n-\pi(n)-2)} \\
     & \sim \frac{2\Big[\frac{3}{{\pi}^2}n^2+n \log \log n-\frac{n^2}{\log n}+\frac{n}{2\log n}\left(\frac{n}{\log n}+1\right)+O(n \log n) \Big]}{\left(n-\frac{n}{\log n}-1\right)\left(n-\frac{n}{\log n}-2\right)}\\
     &= \frac{2\Big[\frac{3}{{\pi}^2}+\frac{\log \log n}{n}-\frac{1}{\log n}+\frac{1}{2\log n}\left(\frac{1}{\log n}+\frac{1}{n}\right)+O(\frac{\log n}{n}) \Big]}
{\left(1-\frac{1}{\log n}-\frac{1}{n}\right)\left(1-\frac{1}{\log n}-\frac{2}{n}\right)}\\
     & = \frac{6}{{\pi}^2}+O\left(\frac{\log n}{n}\right)
  \end{align*}
Therefore, $l(n) \sim \frac{6}{{\pi}^2}$.

\hfill \ensuremath{\Box}

\begin{thm}\label{Theorem 6}
The shortest path length between two nodes in the proposed network with largest possible node index $n$ will be at most $3$ for $n\geq 49$ and it will be at most $2$ for $n\geq 289$.
\end{thm}

\textit{Proof.}
For $n\geq 49$, the graph will be connected from Theorem \ref{Theorem 1}. Now, take two nodes $k=p_1^{a_1}p_2^{a_2}.....p_r^{a_r}$ and $l=p_1^{b_1}p_2^{b_2}.....p_s^{b_s}$, assuming $r\leq s$, with $p_1=2$, $p_2=3$ and so on. If $gcd(k,l)=1$, then they will be connected. If $gcd(k,l)\neq 1$ and $p_{s+1}^2\leq n$ then both $k$ and $l$ will be connected to $p_{s+1}^2$, so the shortest path length will be $2$.

If $gcd(k,l)\neq 1$ and $p_{s+1}^2>n$, then one of the $b_j$'s must be $0$ for $1\leq j\leq s-1$, otherwise it will contradict Theorem \ref{Theorem 1}. If $r=s$, then one of the $a_i$'s will also be $0$ for $1\leq i\leq r-1$ as before, then $k-p_i^2-p_j^2-l$ will be the shortest path from $k$ to $l$ of length $3$ (or of length $2$ if $i=j$) and if $r<s$, then $k-p_t^2-p_j^2-l$ will be the shortest path from $k$ to $l$ of length $3$, where $r+1\leq t\leq s$ and $t\neq j$ (or, again of length $2$ if $t=j$).

Now, for $n\geq 289$, it is enough to show that if $k_n$ and $l_n$ be two nodes as before such that we will always have $p_{s+1}^2\leq n$ when $gcd(k_n,l_n)\neq 1$. Now, in order to maximize $s$, we will consider $k_n=p_1p_2\cdots p_r$ and $l_n=p_1p_{r+1}\cdots p_s$, with $s\geq r+1$, such that $k_n,l_n \leq n$, but $p_1p_2\cdots p_rp_{r+1}>n$ and $p_1p_{r+1}\cdots p_sp_{s+1}>n$, since if we take $l_n=p_1p_t\cdots p_s$, for $t\geq r+2$, then $k_n$ and $l_n$ will be connected to $p_{r+1}^2$. Let for some $r,s\geq 2$, if $n\geq p_{s+1}^2$, then to get an induction argument, if $n+1$ is composite, adding $(n+1)^{th}$ node to the graph, we get two possibilities:

\textbf{Case-1:} Suppose that, $k_{n+1}=k_n=p_1p_2\cdots p_r$, i.e., if $p_1p_2\cdots p_rp_{r+1}>n+1$, then if $p_1p_{r+1}\cdots p_sp_{s+1}>n+1$, we will have $l_{n+1}=l_n=p_1p_{r+1}\cdots p_s$ then from our assumption, we already have $p_{s+1}^2<n+1$, so both $k_{n+1}$ and $l_{n+1}$ will be connected to $p_{s+1}^2$. And if $p_1p_{r+1}\cdots p_sp_{s+1}=l_{n+1}\leq n+1$, then since $p_1p_{r+1}\cdots p_sp_{s+1}>n$ from our assumption, we get $p_1p_{r+1}\cdots p_sp_{s+1}=n+1$. So, using Lemma \ref{Lemma 7} in the appendix, for $s\geq r+2$, we get
\[n+1=p_1p_{r+1}\cdots p_sp_{s+1}\geq p_1p_{s-1}p_sp_{s+1}>p_{s+2}^2.\]
Now, for $s=r+1$,
\[n+1>p_1p_2\cdots p_r=p_1p_2\cdots p_{s-1}>p_{s+2}^2,\]
if $s\geq 6$ from Lemma \ref{Lemma 8} in the Appendix. So $k_{n+1}$ and $l_{n+1}$ both will be connected to $p_{s+2}^2$.

\textbf{Case-2:} In this case, suppose $k_{n+1}\neq k_n$, i.e., if $p_1p_2\cdots p_rp_{r+1}\leq n+1$, then since $p_1p_2\cdots p_rp_{r+1}>n$ from our assumption, we get $k_{n+1}=p_1p_2\cdots p_rp_{r+1}=n+1$. Now, if $l_{n+1}=p_1p_{r+2}\cdots p_s$, then we already have $p_{s+1}^2\leq n<n+1$ from our assumption, so both $k_{n+1}$ and $l_{n+1}$ will be connected to $p_{s+1}^2$. If $l_{n+1}=p_1p_{r+2}\cdots p_sp_{s+1}$, then using Lemma \ref{Lemma 7} in the appendix, for $s\geq r+3$, we get
\[n+1> l_{n+1}=p_1p_{r+2}\cdots p_sp_{s+1}\geq p_1p_{s-1}p_sp_{s+1}>p_{s+2}^2.\]
If $s=r+1$, then
\[n+1=p_1p_2\cdots p_rp_{r+1}=p_1p_2\cdots p_{s-1}p_s>p_1p_2\cdots p_{s-1}>p_{s+2}^2,\]
for $s\geq 6$ from Lemma \ref{Lemma 8} in the appendix (actually $s\geq 5$ is enough for this case). And if $s=r+2$, then
\[n+1=p_1p_2\cdots p_rp_{r+1}=p_1p_2\cdots p_{s-1}>p_{s+2}^2,\]
for $s\geq 6$ from Lemma \ref{Lemma 8} in the appendix, so both $k_{n+1}$ and $l_{n+1}$ will be connected to $p_{s+2}^2$.
And if $l'=p_1p_{r+2}\cdots p_t$, where $t\geq s+2$, then
\[l'=p_1p_{r+2}\cdots p_t\geq p_1p_{r+2}\cdots p_{s+1}p_{s+2}>p_1p_{r+1}\cdots p_{s+1}\geq n+1,\]
since $p_1p_{r+1}.....p_{s+1}>n$ from our assumption. So, $l'$ will not be included to the graph for this case.

So, from both the cases we get $s\geq 6$. Hence, for $n\geq p_{7}^2=289$ the shortest path length will be at most $2$.

\hfill \ensuremath{\Box}

\begin{thm}\label{Theorem 7}
The number of $r$-length labeled cycles, for $r\geq 4$, in the proposed network with largest possible node index $n$ is
$$C_n(r)\leq \left(\frac{6n}{{\pi}^2}\right)^r+O(n^{r-1}\log^{r-1}n).$$
\end{thm}

\textit{Proof.}
We try to count how many $r$-length labeled cycles are there in a coprime network of composite numbers with largest possible node index $n$. We can see that

$$C_n(r) =\sum_{\substack{2\leq a_1,a_2,\cdots ,a_r\leq n \\ a_i\neq ~prime \\ a_i\neq a_j,~i\neq j \\  (a_i,a_{i+1})=1 \\ (a_1,a_r)=1}}\mathbf{1} 
 \leq \sum_{\substack{1\leq a_1,a_2,\cdots ,a_r\leq n \\ (a_i,a_j)=1,\\ i\neq j}}\mathbf{1}$$
 L. T$\acute{o}$th \cite{toth2002probability} showed that
\begin{equation}\label{4}
\sum_{\substack{1\leq a_1,\cdots ,a_r\leq n \\ (a_i,a_j)=1, ~i\neq j \\ (a_i,k)=1}}\mathbf{1}=A_r f_r(k)n^r+O(\theta(k)n^{r-1}\log^{r-1}n),
\end{equation}
where
\[A_r=\prod_{p}\left(1-\frac{1}{p}\right)^{r-1}\left(1+\frac{r-1}{p}\right),\]
\[f_r(k)=\prod_{p|k}\left(1-\frac{r}{p+r-1}\right),\]
and $\theta(k)$ is the number of square-free divisors of $k$.

If the prime factorization of $k$ is $k=p_1^{\alpha_1}p_2^{\alpha_2}\cdots p_{s}^{\alpha_{s}}$, where $s=w(k)$, then
\[\theta(k)=\sum_{i=0}^{s}\binom{s}{i}=2^s=2^{w(k)}.\]
Therefore, using this result with $k=1$ and from Lemma \ref{Lemma 9} in the appendix, we get 
\begin{align*}
\sum_{\substack{1\leq a_1,a_2,\cdots ,a_r\leq n \\ (a_i,a_j)=1,\\ i\neq j}}\mathbf{1} & =n^r\prod_{p}\left(1-\frac{1}{p}\right)^{r-1}\left(1+\frac{r-1}{p}\right) +O(n^{r-1}\log^{r-1}n)\\
& \leq n^r\prod_{p}\left(1-\frac{1}{p^2}\right)^r+O(n^{r-1}\log^{r-1}n)\\
& = n^r\left(\frac{1}{\zeta(2)}\right)^r+O(n^{r-1}\log^{r-1}n)\\
& =\left(\frac{6n}{{\pi}^2}\right)^r+O(n^{r-1}\log^{r-1}n).\\
\end{align*}
Therefore, we get $C_n(r)\leq \left(\frac{6n}{{\pi}^2}\right)^r+O(n^{r-1}\log^{r-1}n).$

\hfill \ensuremath{\Box}

\begin{thm}\label{Theorem 8}
The local clustering coefficient of a node $k$ in the proposed network with largest possible node index $n$ is
\[c(k)\sim \frac{6}{{\pi}^2}\prod_{p|k}\left(\frac{p^2}{p^2-1}\right)\]
\end{thm}

\textit{Proof.}
In order to find the local clustering coefficient of a node $k$, first we have to calculate the total number of triangles $T_n(k)$ which has $k$ as one of the vertices. So,
\[T_n(k)=\frac{1}{2}\sum_{\substack{2\leq i,j\leq n \\ i,j \neq prime \\ (i,k)=1 \\ (j,k)=1 \\ (i,j)=1}}\mathbf{1}.\]
Now,
\begin{equation}\label{}
\sum_{\substack{2\leq i,j\leq n \\ i,j \neq prime \\ (i,k)=1 \\ (j,k)=1 \\ (i,j)=1}}\mathbf{1}
=\sum_{\substack{1\leq i,j\leq n \\ (i,k)=1 \\ (j,k)=1 \\ (i,j)=1}}\mathbf{1}
-2\sum_{\substack{i=1\\ (i,k)=1}}^{n}\sum_{\substack{p\leq n\\ p\nmid i\\ p\nmid k}}\mathbf{1}
+\sum_{\substack{p\leq n\\  p\nmid k}}\sum_{\substack{q\leq n\\ q\nmid i\\ q\neq p}}\mathbf{1}
-2\sum_{\substack{i=1\\ (i,k)=1}}^{n}\mathbf{1}
+2\sum_{\substack{p\leq n\\ p\nmid k}}\mathbf{1}+1.
\end{equation}
We can find the closed expressions for all the sums easily except the first one, say $S_1$. So we get
\begin{equation}\label{}
\sum_{\substack{2\leq i,j\leq n \\ i,j \neq prime \\ (i,k)=1 \\ (j,k)=1 \\ (i,j)=1}}\mathbf{1}
=\sum_{\substack{1\leq i,j\leq n \\ (i,k)=1 \\ (j,k)=1 \\ (i,j)=1}}\mathbf{1}
+2\sum_{\substack{i=1 \\ (i,k)=1}}^{n}w(i)
-\left[\pi(n)-w(k)\right]\left[2\phi(n,k)-\pi(n)+w(k)-1\right]-2\phi(n,k)+1.
\end{equation}

We can see that $S_1$ is the $r=2$ case of (\ref{4}). Since,
\[A_2=\prod_{p}\left(1-\frac{1}{p}\right)\left(1+\frac{1}{p}\right)=\prod_{p}\left(1-\frac{1}{p^2}\right)=\frac{1}{\zeta(2)}=\frac{6}{{\pi}^2},\]
\[f_2(k)=\prod_{p|k}\left(1-\frac{2}{p+1}\right)=\prod_{p|k}\left(\frac{p-1}{p+1}\right),\]
So, we get
\[S_1=\frac{6}{{\pi}^2}\prod_{p|k}\left(\frac{p-1}{p+1}\right)n^2+O(2^{w(k)}n\log{n}).\]
Therefore,
\begin{equation}\label{}
\begin{split}
   T_n(k)= & \frac{1}{2}\Big[\frac{6}{{\pi}^2}\prod_{p|k}\left(\frac{p-1}{p+1}\right)n^2-\left[\pi(n)-w(k)\right]\left[2\phi(n,k)-\pi(n)+w(k)-1\right]\\
     & ~~~~~~~~~~~~~~~~~+2\sum_{\substack{i=1 \\ (i,k)=1}}^{n}w(i)-2\phi(n,k)+O(2^{w(k)}n\log{n})\Big].
\end{split}
\end{equation}
So, the local clustering coefficient of node $k$,
\begin{equation}\label{}
 c(k)=\frac{T_n(k)}{\binom{d_k}{2}}.
\end{equation}
Now,
\[T_n(k)\sim \frac{1}{2}\Big[\frac{6}{{\pi}^2}\prod_{p|k}\left(\frac{p-1}{p+1}\right)n^2-\left[\frac{n}{\log n}-w(k)\right]\left[\frac{2n\phi(k)}{k}-\frac{n}{\log n}+w(k)-1\right]+O(2^{w(k)}n\log{n})\Big],\]
and
\begin{align*}
  \binom{d_k}{2} & =\frac{(\phi(n,k)-\pi(n)+w(k)-1)(\phi(n,k)-\pi(n)+w(k)-2)}{2} \\
   & \sim \frac{(\frac{n\phi(k)}{k}-\frac{n}{\log n}+w(k)-1)(\frac{n\phi(k)}{k}-\frac{n}{\log n}+w(k)-2)}{2}.
\end{align*}

Therefore,
\begin{align*}
  c(k)& \sim \frac{\Big[\frac{6}{{\pi}^2}\prod_{p|k}\left(\frac{p-1}{p+1}\right)n^2-\left[\frac{n}{\log n}-w(k)\right]\left[\frac{2n\phi(k)}{k}-\frac{n}{\log n}+w(k)-1\right]+O(2^{w(k)}n\log{n})\Big]}{(\frac{n\phi(k)}{k}-\frac{n}{\log n}+w(k)-1)(\frac{n\phi(k)}{k}-\frac{n}{\log n}+w(k)-2)} \\
   & = \frac{\Big[\frac{6}{{\pi}^2}\prod_{p|k}\left(\frac{p-1}{p+1}\right)-\left[\frac{1}{\log n}-\frac{w(k)}{n}\right]\left[\frac{2\phi(k)}{k}-\frac{1}{\log n}+\frac{w(k)}{n}-\frac{1}{n}\right]+O(2^{w(k)}\frac{\log{n}}{n})\Big]}{(\frac{\phi(k)}{k}-\frac{1}{\log n}+\frac{w(k)}{n}-\frac{1}{n})(\frac{\phi(k)}{k}-\frac{1}{\log n}+\frac{w(k)}{n}-\frac{2}{n})} \\
   & \sim \frac{6}{{\pi}^2}\prod_{p|k}\left(\frac{p-1}{p+1}\right)\left(\frac{k}{\phi(k)}\right)^2.
\end{align*}
Since $\phi(k)=k\prod_{p|k}(1-\frac{1}{p})$, so finally we have the form of $c(k)$ as
\begin{equation}\label{}
 c(k)\sim \frac{6}{{\pi}^2}\prod_{p|k}\left(\frac{p^2}{p^2-1}\right).
\end{equation}

\hfill \ensuremath{\Box}

\begin{figure}[htbp]
    \centering
        \includegraphics[scale=.6]{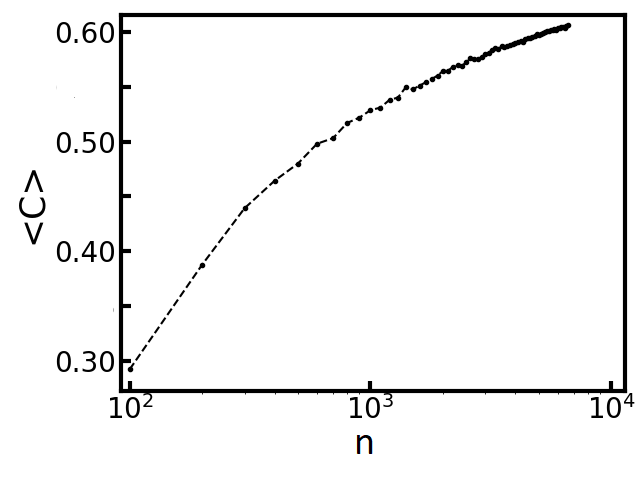}
       \caption{{The average of the local clustering coefficients of the nodes is plotted by changing the highest possible node index $n$. It shows the average of the local clustering coefficients approaches near $\frac{6}{{\pi}^2}(\approx 0.61)$ with the increase in $n$. We have used log-scale for $x$-axis.}}
    \label{Fig:3}
    \end{figure}

\section{Weak Pseudo-randomness}
In this section, we show that the proposed network is weakly pseudo-random. Apart from the actual definition of weakly pseudo-random graphs, that is for a real parameter $p = p(n)\in (0,1)$, a sequence of graphs $(G_n)$ is weakly pseudo-random \cite{krivelevich2006pseudo} if for all subsets $U\subseteq V(G_n)$, we have
$$\left|E(U)-p\binom{|U|}{2}\right|=o(pn^2),$$
we can also find some equivalent conditions to this in literature \cite{chung1989quasi}. In particular, we are interested in the following statement \cite{krivelevich2006pseudo}: A graph is weakly pseudo-random if and only if 
$$\sum_{x,y\in \mathbf{V}}\left|codeg(x,y)-p^2N\right|=o(N^3),$$
where $N$ is the number of nodes and  $codeg(x,y)$ is the number of common neighbors of nodes $x$ and $y$.

Now, since 
$$codeg(x,y)=\sum_{i\in \mathbf{V}}\mathbf{1}[x\sim i]\mathbf{1}[y\sim i]=\sum_{i\in \mathbf{V}}A[x,i]A[y,i]=A^2[x,y],$$
where $A$ is the adjacency matrix of the graph. So, using Theorem \ref{Theorem 7}, we get
$$\sum_{x,y\in \mathbf{V}}codeg(x,y)^2=\sum_{x,y\in \mathbf{V}}[A^2\left[x,y]\right]^2=Tr(A^4)=C_n(4)\leq \left(\frac{6n}{{\pi}^2}\right)^4+O(n^3\log^3 n)=(pn)^4+o(n^4),$$
where $p=\frac{6}{{\pi}^2}$. On the other hand,
\begin{align*}
 \sum_{x,y\in \mathbf{V}}codeg(x,y) &  =\sum_{x,y\in \mathbf{V}}\sum_{i\in \mathbf{V}}\mathbf{1}[x\sim i]\mathbf{1}[y\sim i]=\sum_{i\in \mathbf{V}}\left[\sum_{x\in \mathbf{V}}\mathbf{1}[x\sim i]\right]^2\\
 &=\sum_{i\in \mathbf{V}}deg(i)^2\geq \frac{\left[\sum_{i\in \mathbf{V}}deg(i)\right]^2}{N(n)}=N(n)\overline{d(n)}^2\\
 &=(n-\log n-1)\left[\frac{6}{{\pi}^2}n+O\left(\frac{n}{(\log n)^2}\right)\right]^2=p^2n^3+o(n^3).
\end{align*}
Now, using Cauchy-Schwarz inequality, we get

$$\sum_{x,y\in \mathbf{V}}\left|codeg(x,y)-p^2N(n)\right| \leq N(n)\left[\sum_{x,y\in \mathbf{V}}\left(codeg(x,y)-p^2N(n)\right)^2\right]^{1/2}$$
$$= N(n)\left[ \sum_{x,y\in \mathbf{V}}codeg(x,y)^2-2p^2N(n)\sum_{x,y\in \mathbf{V}}codeg(x,y)+p^4N(n)^4\right]^{1/2}$$
$$\leq N(n)\left[(pn)^4+o(n^4)-2p^2N(n)(p^2n^3+o(n^3))+p^4N(n)^4\right]^{1/2}.$$
Therefore,
$$\sum_{x,y\in \mathbf{V}}\left|codeg(x,y)-p^2N(n)\right|=o(N(n)^3).$$
Therefore, this sequence of graphs is weakly pseudo-random with $p=\frac{6}{{\pi}^2}$.

Now, we can directly use the other equivalent conditions\cite{krivelevich2006pseudo} to the weak pseudo-randomness and infer other properties of this network. For instance, given a graph $L$ of $l\geq 4$ vertices, let $N_G^{*}(L)$ be the number of labeled induced copies of $L$ in a weakly pseudo-random graph, $G$ with $n$ vertices, then
\begin{equation}\label{5}
N_G^{*}(L)=(1+o(1))n^lp^{|E(L)|}(1-p)^{\binom{l}{2}-|E(L)|}.
\end{equation}
Taking $L$ to be $r$-length cycle for $r\geq 4$ in (\ref{5}), we get
$$C_n(r)=(1+o(1))(N(n))^r p^{r}(1-p)^{\frac{r(r-3)}{2}}.$$
So, we get a cycle of length $r$ if $N(n)p(1-p)^{\frac{(r-3)}{2}}\geq 1$. This indicates the existence of a cycle of length $O(\log n)$ for large $n$ as $p=\frac{6}{{\pi}^2}$ is constant.

\begin{figure}[htbp]
    \centering
        \includegraphics[scale=.6]{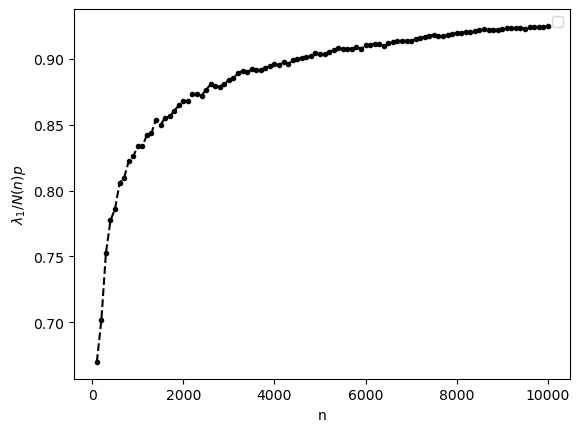}
       \caption{{The ratio of the maximum eigenvalue $\lambda_1$ of the adjacency matrix and N(n)p is plotted by changing the highest possible node index $n$. The ratio approaching 1 provides evidence to $\lambda_1=(1+o(1))N(n)p$.}}
    \label{Fig:6}
    \end{figure}

Again from \cite{krivelevich2006pseudo}, we get weak pseudo-randomness implies that the maximum eigenvalue of the adjacency matrix, $\lambda_1=(1+o(1))N(n)p$, which is also evident from Fig:\ref{Fig:6}.

\section{Discussion and possible applications}

In this paper, we have constructed a network with nodes indexed by the composite numbers and we defined the adjacency of two nodes when their indices are coprime to each other. As the size of the network increases by adding nodes one by one to the network starting from the one with index 4, we have noticed a number of interesting behaviors in some of the statistics that we have found. We have seen that the link density saturates to $\frac{6}{{\pi}^2}$ and the average degree of the nodes increase linearly with the increment in $n$, which is the highest possible node index, with slope again $\frac{6}{{\pi}^2}$. This behavior actually supports a well known result in Number Theory that the probability that two randomly chosen natural numbers are relatively prime is $\frac{6}{{\pi}^2}$ \cite{apostol1998introduction}. It also refers that for large $n$, $\frac{6}{{\pi}^2}$(i.e., $60.79\%$) of the total number of possible edges will be present in the network. We have also seen that the network will not have any isolated vertex when $n\geq49$ that means after this point the network will be connected and starting from any vertex it will always be possible to visit all the other vertices in this network. But starting from a given vertex, what is the minimum number of vertices we need to visit in order to reach a target vertex? We have proved in this context that the shortest path length between any two nodes will be at most 3 when $49\leq n\leq288$ and it is at most 2 for $n\geq 289$. A very interesting graph property is the existence of Hamiltonian cycles of a given length, that is the cycles of a specific length that doesn't visit any vertex more than once. We have found that the number of these cycles of $r$-length is $(1+o(1))(N(n))^r p^{r}(1-p)^{\frac{r(r-3)}{2}}$ for large $n$. But for a fixed $n$, as we increase the length $r$ of the cycles, this asymptotic number of cycles starts dropping and eventually becomes less than $1$, which suggests us to investigate what is the maximum possible length of a cycle that may exist in the proposed network. We found that there exists a cycle of length at most $O(\log n)$ for large $n$. Also for a given vertex, while investigating the connectivity between its neighbors, we have found a measure named local clustering coefficient $c(k)$ of a node $k$, which is the proportion of the total number of triangle whose one vertex is that node itself to the total number of possible triangles with that node as a vertex. We have shown that for large $n$, 

$$c(k)\sim \frac{6}{{\pi}^2}\prod_{p|k}\left(\frac{p^2}{p^2-1}\right),$$
which suggests as the size of the network increases, more than $\frac{6}{{\pi}^2} (i.e.,~60.79\%$) of the neighbors of any node will be connected to each other.

Though the construction of this network being completely deterministic, it's study will still be as interesting and as important as the other random networks because of the irregular distribution of primes over natural numbers which makes the coprimality of two nodes to be irregular too. Following Krivelevich and Sudakov \cite{krivelevich2006pseudo} and using some previously discussed results, we have shown that the graph sequence namely coprime networks of composite numbers is weakly pseudo-random, which again supports the property that two integers being relative prime is random-like and that is actually a very strong motif in analytic number theory. 

\begin{figure}[htbp]
    \centering
    \includegraphics[scale=.77]{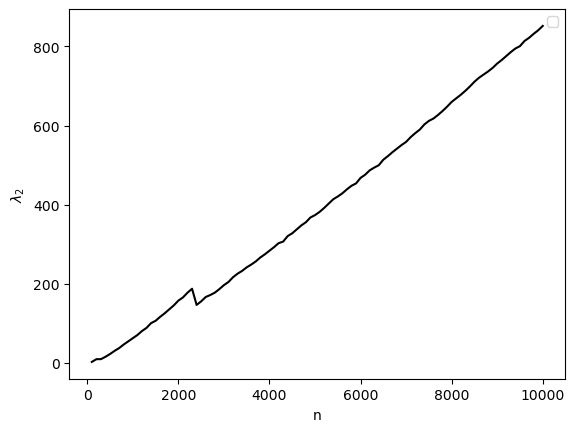}
        \caption{{The variation of the second eigenvalue $\lambda_2$ of the Laplacian matrix is plotted by changing the highest possible node index $n$. It is observed that the value of $\lambda_2$ is increasing with the increased value of node indices, which provides evidence of the high connectivity in the constructed network.}}
    \label{eigenvalues}
    \end{figure}
\begin{figure}[htbp]
    \centering
        \includegraphics[scale=.77]{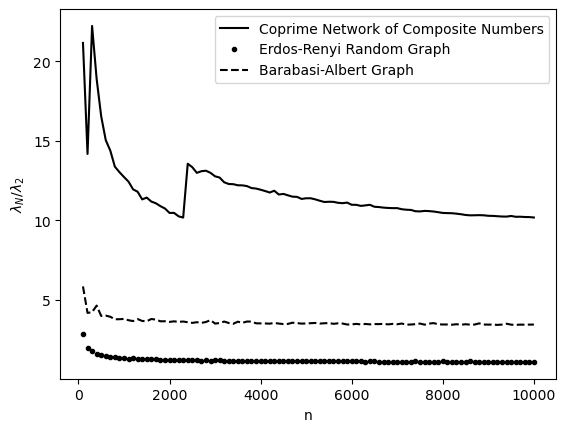}
        \caption{{The ratio of $N$th and 2nd eigenvalues of the Laplacian matrix is plotted with respect to the highest possible node index $n$ for the constructed coprime network of composite numbers and is compared with Erd\H{o}s–R\'{e}nyi random network and Barab\'{a}si-Albert scale-free network with the same number of nodes $N(n)$ and edges $E(n)$. This figure suggests less synchronizability in the proposed network as compared to the other two.}}
    \label{eigenvalues ratio}
    \end{figure}
One of the possible application to use any network configuration is to study collective behavior, namely synchronization, in dynamical networks. The synchronizability of a network is depended on the eigenvalues of the zero-row-sum real Laplacian matrix $L$. If $L$ is a $N \times N$ matrix then its eigenvalues satisfies $0=\lambda_1\le \lambda_2\le \cdots \le \lambda_N.$ The smaller values of $\frac{\lambda_N}{\lambda_2}$ gives the more synchronizability in the network and vice-versa \cite{jafarizadeh2020optimizing, tang2014synchronization,  barahona2002synchronization}. We first calculate the second eigenvalue $\lambda_2$ of $L$ of the constructed coprime network of composite numbers by changing the highest possible node index $n$, shown in Fig.\ref{eigenvalues}. From this figure, it is noticed that the value of $\lambda_2$ increases as we increase the number of nodes and the corresponding links, which gives an evidence of the high connectivity in the constructed network. Then we calculate the ratio $\frac{\lambda_N}{\lambda_2}$ in order to see the synchronizability in the network by varying the highest possible node index $n$. Next we compare the synchronizability of the constructed coprime network of composite numbers with two other networks, namely Erd\H{o}s–R\'{e}nyi random network and Barab\'{a}si-Albert network with the same number of nodes $N(n)$ and edges $E(n)$. Comparing with these two networks, we can notice that the synchronizability in the coprime networks of composite numbers is less from Fig.\ref{eigenvalues ratio}. This less synchronizability is very much desirable in ecological and predator-prey networks to increase prolonged transient dynamics \cite{holland2008strong}.

\section{Conclusion}
In this paper, we have studied a network which we have constructed by taking all the composite numbers in a given range starting from 4 and drawing an undirected link between two of them if they are coprime to each other. Therefore, the nodes having same set of prime divisors will have exactly same neighbors, hence same degree and same value of local clustering coefficients. Among them, the prime powers have larger number of neighbors. We have shown that the network will be connected if the largest possible node index $n\geq 49$ and the maximum of the shortest path lengths is 3 if $n$ lies between 49 and 288, whereas it is 2 if $n\geq 289$. We have also derived the analytic expressions for the maximum degree, the average degrees of the nodes, link density and local clustering coefficients of the nodes. Also, we have found an expression for the number of $r$-length labeled cycles. Finally, we show that the constructed network is weakly pseudo-random. The root cause of this is the coprimality of two integers is random-like.

Also, we have seen that this network has less synchronizability which is noticed in ecological and predator-prey
networks. Someone can also try to find the degree distribution, i.e., if we choose a node randomly, then what is the probability that it would have degree $m$,
$$\mathbb{P}(m) =\frac{|\{k\in V:d_k=m\}|}{n-\pi(n)-1}
=\frac{|\{k\in V:\phi(n,k)+w(k)=m+\pi(n)+1\}|}{n-\pi(n)-1}.$$

\section*{Declaration of competing interest}
The authors declare that they have no known competing financial interests or personal relationships that could have appeared to influence the work reported in this paper.

\section*{Data availability}
No data was used for the research described in the article.

\section*{Remarks}
The codes of the plots can be found in \href{https://github.com/Rahil-Miraj/Coprime-Networks-of-Composite-Numbers.git}{https://github.com/Rahil-Miraj/Coprime-Networks-of-Composite-Numbers.git}.

\section*{Acknowledgments}
The authors would like to thank Sudipta Mallik, Parthanil Roy, B. Sury and  Niranjan Balachandran for their suggestions and Subrata Ghosh for helping in the simulations.

\section*{Appendix}

In this appendix, we state some results which we have used in the context of proving the results throughout the paper.

\begin{lem}\label{Lemma 1}
For $x>1$, we have \cite{apostol1998introduction}
\[\sum_{k\leq x}\phi(k)=\frac{3}{{\pi}^2}x^2+O(x \log x).\]
\end{lem}

\begin{lem}\label{Lemma 2}
For $x\geq 1$, we have \cite{hardy1979introduction}
\[\sum_{k\leq x}\omega(k)=x \log \log x+B_1x+O(x),\]
where  $B_{1}\approx 0.2614972128$ is the Mertens constant.
\end{lem}

\begin{lem} \label{Lemma 3}
For $x\geq 1$, we have \cite{hardy1979introduction}
\[\sum_{k\leq x}\omega(k)^2=x (\log \log x)^2+O(x\log \log x).\]
\end{lem}
\begin{lem} \label{Lemma 5}
For $x>2$, we have
\[\sum_{k\leq x}\pi(k)={\lfloor x \rfloor} \pi(x) -\sum_{p\leq x}p.\]
\end{lem}

\textit{Proof.}
Let $m=\pi(x)$. Then,
\begin{align*}
\sum_{k\leq x}\pi(k) &=\sum_{k=2}^{{\lfloor x \rfloor}}\pi(k)=\sum_{i=2}^{m}(p_i-p_{i-1})\pi(p_i)+({\lfloor x \rfloor}-p_m)(\pi(p_m)+1)\\
&=\sum_{i=2}^{m}(p_i-p_{i-1})(i-1)+({\lfloor x \rfloor}-p_m)((m-1)+1)\\
&=\sum_{i=2}^{m}ip_i-\sum_{i=2}^{m}(i-1)p_{i-1}-\sum_{i=2}^{m}p_i+m({\lfloor x \rfloor}-p_m)\\
&=m{\lfloor x \rfloor}-\sum_{i=1}^{m}p_i.
\end{align*}

\hfill \ensuremath{\Box}

\begin{lem}\label{Lemma 6}
Let $t\geq 2$ be an integer. Then, for $s\geq 2$, whenever $p_1p_tp_{t+1}.....p_s> p_{s+1}^2$ with $p_1=2, p_2=3$ and so on, we will have $p_1p_tp_{t+1}.....p_{s+1}> p_{s+2}^2$.
\end{lem}

\textit{Proof.}
Let $p_1p_tp_{t+1}.....p_s> p_{s+1}^2$. Then,
\begin{equation}\label{1}
p_1p_t.....p_sp_{s+1}> p_{s+1}^3.
\end{equation}

Now, from Bertrand's postulate \cite{sondow2009ramanujan}, there is a prime between $x$ and $2x$ for all $x>1$, so
\[p_{s+1}<p_{s+2}<2p_{s+1},\]
\[p_{s+2}^2<4p_{s+1}^2<p_{s+1}^3,\]
if $p_{s+1}>4$, i.e., $p_{s+1}\geq 5=p_3$.

Therefore, (\ref{1}) gives
\[p_1p_tp_{t+1}.....p_{s+1}> p_{s+2}^2,\]
for $s\geq 2$. 

\hfill \ensuremath{\Box}

\begin{lem}\label{Lemma 7}
 Let $t\geq 3$ be an integer. Then, $p_1p_{t-1}p_tp_{t+1}>p_{t+2}^2$.
\end{lem}

\textit{Proof.}
We can easily check that this holds true for $t=3$ and $t=4$. Now, using Bertrand's postulate again, we get $p_1p_{t+1}=2p_{t+1}>p_{t+2}$. Now, since $11=p_5$ is the $2^{nd}$ Ramanujan prime 
 \cite{sondow2009ramanujan}, $p_{t-1}p_t>2p_t>p_{t+2}$ for $t\geq 5$. So, combining these two results, we get $p_1p_{t-1}p_tp_{t+1}>p_{t+2}^2$ for $t\geq 5$.

\hfill \ensuremath{\Box}

\begin{lem}\label{Lemma 8}
 For an integer $t\geq 6$, $p_1p_2....p_{t-1}>p_{t+2}^2$.
\end{lem}

\textit{Proof.}
Let us consider the sequence, $\{a_t\}_{t=2}^{\infty}$, defined as $a_t=p_1p_2....p_{t-1}-p_{t+2}^2$. Then we will have
\[a_{t+1}=p_1p_2....p_t-p_{t+3}^2>p_1p_2....p_t-p_tp_{t+2}^2=p_t(p_1p_2....p_{t-1}-p_{t+2}^2)=p_ta_t>a_t,\]
for $t>3$, since, $p_tp_{t+2}^2>4p_{t+2}^2>p_{t+3}^2$, again from Bertrand's postulate. Therefore, $\{a_t\}$ is increasing for $t\geq 4$ and since we get the first positive value of $a_t$ at $t=6$, we conclude that for $t\geq 6$, we will always get $p_1p_2....p_{t-1}>p_{t+2}^2$.

\hfill \ensuremath{\Box}

\begin{lem}\label{Lemma 9}
 For $x\geq 2$ and $r\geq 4$, we have
 $$\left(1-\frac{1}{x}\right)^{r-1}\left(1+\frac{r-1}{x}\right)\leq \left(1-\frac{1}{x^2}\right)^r.$$
\end{lem}

\textit{Proof.}
It is enough to show that 
$$\left(1+\frac{r-1}{x}\right)\leq \left(1-\frac{1}{x}\right)\left(1+\frac{1}{x}\right)^r$$
holds for $x\geq 2$ and $r\geq 4$.

Let $f(r,x)=\left(1-\frac{1}{x}\right)\left(1+\frac{1}{x}\right)^r-\left(1+\frac{r-1}{x}\right)$. We will show that  $f(r,x)$ is increasing in $r$ for $r\geq 4$ and $x\geq 2$ and $\min_{r\geq 4}f(r,x)=f(4,x)\geq 0$ for $x\geq 2$.

We have 
$$\frac{\partial f}{\partial r}=\left(1-\frac{1}{x}\right)\left(1+\frac{1}{x}\right)^r\log r-1/x>1-1/x>0,$$
for $r\geq 4$ and $x\geq 2$. Hence, $f(r,x)$ is increasing in $r$ for $r\geq 4$ and $x\geq 2$.

Now, suppose for some $x\geq 2$, we will get $f(4,x)<0$. Then since $f(4,2)=0.0313>0$, $f(4,x)=\frac{1}{x^5}(2x^3-2x^2-3x-1)$ will have a zero in $x>2$. Then $2x^3-2x^2-3x-1$ will have a zero in $x>2$. But $2x^3-2x^2-3x-1$ is increasing for $x>2$ as $\frac{d(2x^3-2x^2-3x-1)}{dx}=6x^2-4x-3>0$ and takes value $1$ for $x=2$, a contradiction.

\hfill \ensuremath{\Box}

\bibliographystyle{elsarticle-num-names}

\begin{thebibliography}{47}
	\expandafter\ifx\csname natexlab\endcsname\relax\def\natexlab#1{#1}\fi
	\providecommand{\url}[1]{\texttt{#1}}
	\providecommand{\href}[2]{#2}
	\providecommand{\path}[1]{#1}
	\providecommand{\DOIprefix}{doi:}
	\providecommand{\ArXivprefix}{arXiv:}
	\providecommand{\URLprefix}{URL: }
	\providecommand{\Pubmedprefix}{pmid:}
	\providecommand{\doi}[1]{\href{http://dx.doi.org/#1}{\path{#1}}}
	\providecommand{\Pubmed}[1]{\href{pmid:#1}{\path{#1}}}
	\providecommand{\bibinfo}[2]{#2}
	\ifx\xfnm\relax \def\xfnm[#1]{\unskip,\space#1}\fi
	\bibitem[{Newman(2010)}]{Newman_Network_Book}
	\bibinfo{author}{M.~Newman}, \bibinfo{title}{Networks: an introduction},
	\bibinfo{publisher}{Oxford University Press}, \bibinfo{year}{2010}.
	\bibitem[{Albert and Barab{\'a}si(2002)}]{albert2002statistical}
	\bibinfo{author}{R.~Albert}, \bibinfo{author}{A.-L. Barab{\'a}si},
	\newblock \bibinfo{title}{Statistical mechanics of complex networks},
	\newblock \bibinfo{journal}{Reviews of Modern Physics} \bibinfo{volume}{74}
	(\bibinfo{year}{2002}) \bibinfo{pages}{47}.
	\bibitem[{Pastor-Satorras et~al.(2003)Pastor-Satorras, Rubi, and
		Diaz-Guilera}]{pastor2003statistical}
	\bibinfo{author}{R.~Pastor-Satorras}, \bibinfo{author}{M.~Rubi},
	\bibinfo{author}{A.~Diaz-Guilera}, \bibinfo{title}{Statistical mechanics of
		complex networks}, volume \bibinfo{volume}{625}, \bibinfo{publisher}{Springer
		Science \& Business Media}, \bibinfo{year}{2003}.
	\bibitem[{Boccaletti et~al.(2006)Boccaletti, Latora, Moreno, Chavez, and
		Hwang}]{boccaletti2006complex}
	\bibinfo{author}{S.~Boccaletti}, \bibinfo{author}{V.~Latora},
	\bibinfo{author}{Y.~Moreno}, \bibinfo{author}{M.~Chavez},
	\bibinfo{author}{D.-U. Hwang},
	\newblock \bibinfo{title}{Complex networks: Structure and dynamics},
	\newblock \bibinfo{journal}{Physics Reports} \bibinfo{volume}{424}
	(\bibinfo{year}{2006}) \bibinfo{pages}{175--308}.
	\bibitem[{Arenas et~al.(2008)Arenas, D{\'\i}az-Guilera, Kurths, Moreno, and
		Zhou}]{arenas2008synchronization}
	\bibinfo{author}{A.~Arenas}, \bibinfo{author}{A.~D{\'\i}az-Guilera},
	\bibinfo{author}{J.~Kurths}, \bibinfo{author}{Y.~Moreno},
	\bibinfo{author}{C.~Zhou},
	\newblock \bibinfo{title}{Synchronization in complex networks},
	\newblock \bibinfo{journal}{Physics Reports} \bibinfo{volume}{469}
	(\bibinfo{year}{2008}) \bibinfo{pages}{93--153}.
	\bibitem[{Barrat et~al.(2008)Barrat, Barthelemy, and
		Vespignani}]{barrat2008dynamical}
	\bibinfo{author}{A.~Barrat}, \bibinfo{author}{M.~Barthelemy},
	\bibinfo{author}{A.~Vespignani}, \bibinfo{title}{Dynamical processes on
		complex networks}, \bibinfo{publisher}{Cambridge University Press},
	\bibinfo{year}{2008}.
	\bibitem[{Cohen and Havlin(2010)}]{Cohen_complex_book}
	\bibinfo{author}{R.~Cohen}, \bibinfo{author}{S.~Havlin},
	\bibinfo{title}{Complex networks: structure, robustness and function},
	\bibinfo{publisher}{Cambridge University Press}, \bibinfo{year}{2010}.
	\bibitem[{Battiston et~al.(2021)Battiston, Amico, Barrat, Bianconi, Ferraz~de
		Arruda, Franceschiello, Iacopini, K{\'e}fi, Latora, Moreno
		et~al.}]{battiston2021physics}
	\bibinfo{author}{F.~Battiston}, \bibinfo{author}{E.~Amico},
	\bibinfo{author}{A.~Barrat}, \bibinfo{author}{G.~Bianconi},
	\bibinfo{author}{G.~Ferraz~de Arruda}, \bibinfo{author}{B.~Franceschiello},
	\bibinfo{author}{I.~Iacopini}, \bibinfo{author}{S.~K{\'e}fi},
	\bibinfo{author}{V.~Latora}, \bibinfo{author}{Y.~Moreno}, et~al.,
	\newblock \bibinfo{title}{The physics of higher-order interactions in complex
		systems},
	\newblock \bibinfo{journal}{Nature Physics} \bibinfo{volume}{17}
	(\bibinfo{year}{2021}) \bibinfo{pages}{1093--1098}.
	\bibitem[{Hens et~al.(2019)Hens, Harush, Haber, Cohen, and
		Barzel}]{hens2019spatiotemporal}
	\bibinfo{author}{C.~Hens}, \bibinfo{author}{U.~Harush},
	\bibinfo{author}{S.~Haber}, \bibinfo{author}{R.~Cohen},
	\bibinfo{author}{B.~Barzel},
	\newblock \bibinfo{title}{Spatiotemporal signal propagation in complex
		networks},
	\newblock \bibinfo{journal}{Nature Physics} \bibinfo{volume}{15}
	(\bibinfo{year}{2019}) \bibinfo{pages}{403--412}.
	\bibitem[{Gao et~al.(2016)Gao, Barzel, and Barab{\'a}si}]{gao2016universal}
	\bibinfo{author}{J.~Gao}, \bibinfo{author}{B.~Barzel}, \bibinfo{author}{A.-L.
		Barab{\'a}si},
	\newblock \bibinfo{title}{Universal resilience patterns in complex networks},
	\newblock \bibinfo{journal}{Nature} \bibinfo{volume}{530}
	(\bibinfo{year}{2016}) \bibinfo{pages}{307--312}.
	\bibitem[{Barzel and Barab{\'a}si(2013)}]{barzel2013universality}
	\bibinfo{author}{B.~Barzel}, \bibinfo{author}{A.-L. Barab{\'a}si},
	\newblock \bibinfo{title}{Universality in network dynamics},
	\newblock \bibinfo{journal}{Nature Physics} \bibinfo{volume}{9}
	(\bibinfo{year}{2013}) \bibinfo{pages}{673--681}.
	\bibitem[{Meena et~al.(2023)Meena, Hens, Acharyya, Haber, Boccaletti, and
		Barzel}]{meena2023emergent}
	\bibinfo{author}{C.~Meena}, \bibinfo{author}{C.~Hens},
	\bibinfo{author}{S.~Acharyya}, \bibinfo{author}{S.~Haber},
	\bibinfo{author}{S.~Boccaletti}, \bibinfo{author}{B.~Barzel},
	\newblock \bibinfo{title}{Emergent stability in complex network dynamics},
	\newblock \bibinfo{journal}{Nature Physics} \bibinfo{volume}{19}
	(\bibinfo{year}{2023}) \bibinfo{pages}{1033--1042}.
	\bibitem[{Artime et~al.(2024)Artime, Grassia, De~Domenico, Gleeson, Makse,
		Mangioni, Perc, and Radicchi}]{artime2024robustness}
	\bibinfo{author}{O.~Artime}, \bibinfo{author}{M.~Grassia},
	\bibinfo{author}{M.~De~Domenico}, \bibinfo{author}{J.~P. Gleeson},
	\bibinfo{author}{H.~A. Makse}, \bibinfo{author}{G.~Mangioni},
	\bibinfo{author}{M.~Perc}, \bibinfo{author}{F.~Radicchi},
	\newblock \bibinfo{title}{Robustness and resilience of complex networks},
	\newblock \bibinfo{journal}{Nature Reviews Physics}  (\bibinfo{year}{2024})
	\bibinfo{pages}{1--18}.
	\bibitem[{Borgatti et~al.(2009)Borgatti, Mehra, Brass, and
		Labianca}]{borgatti2009network}
	\bibinfo{author}{S.~P. Borgatti}, \bibinfo{author}{A.~Mehra},
	\bibinfo{author}{D.~J. Brass}, \bibinfo{author}{G.~Labianca},
	\newblock \bibinfo{title}{Network analysis in the social sciences},
	\newblock \bibinfo{journal}{Science} \bibinfo{volume}{323}
	(\bibinfo{year}{2009}) \bibinfo{pages}{892--895}.
	\bibitem[{Menck et~al.(2013)Menck, Heitzig, Marwan, and
		Kurths}]{menck2013basin}
	\bibinfo{author}{P.~J. Menck}, \bibinfo{author}{J.~Heitzig},
	\bibinfo{author}{N.~Marwan}, \bibinfo{author}{J.~Kurths},
	\newblock \bibinfo{title}{How basin stability complements the linear-stability
		paradigm},
	\newblock \bibinfo{journal}{Nature Physics} \bibinfo{volume}{9}
	(\bibinfo{year}{2013}) \bibinfo{pages}{89--92}.
	\bibitem[{Wang et~al.(2019)Wang, Liu, Liang, Hu, and
		Zhou}]{wang2019coevolution}
	\bibinfo{author}{W.~Wang}, \bibinfo{author}{Q.-H. Liu},
	\bibinfo{author}{J.~Liang}, \bibinfo{author}{Y.~Hu},
	\bibinfo{author}{T.~Zhou},
	\newblock \bibinfo{title}{Coevolution spreading in complex networks},
	\newblock \bibinfo{journal}{Physics Reports} \bibinfo{volume}{820}
	(\bibinfo{year}{2019}) \bibinfo{pages}{1--51}.
	\bibitem[{Pastor-Satorras et~al.(2015)Pastor-Satorras, Castellano, Van~Mieghem,
		and Vespignani}]{pastor2015epidemic}
	\bibinfo{author}{R.~Pastor-Satorras}, \bibinfo{author}{C.~Castellano},
	\bibinfo{author}{P.~Van~Mieghem}, \bibinfo{author}{A.~Vespignani},
	\newblock \bibinfo{title}{Epidemic processes in complex networks},
	\newblock \bibinfo{journal}{Reviews of Modern Physics} \bibinfo{volume}{87}
	(\bibinfo{year}{2015}) \bibinfo{pages}{925}.
	\bibitem[{Donohue et~al.(2016)Donohue, Hillebrand, Montoya, Petchey, Pimm,
		Fowler, Healy, Jackson, Lurgi, McClean et~al.}]{donohue2016navigating}
	\bibinfo{author}{I.~Donohue}, \bibinfo{author}{H.~Hillebrand},
	\bibinfo{author}{J.~M. Montoya}, \bibinfo{author}{O.~L. Petchey},
	\bibinfo{author}{S.~L. Pimm}, \bibinfo{author}{M.~S. Fowler},
	\bibinfo{author}{K.~Healy}, \bibinfo{author}{A.~L. Jackson},
	\bibinfo{author}{M.~Lurgi}, \bibinfo{author}{D.~McClean}, et~al.,
	\newblock \bibinfo{title}{Navigating the complexity of ecological stability},
	\newblock \bibinfo{journal}{Ecology Letters} \bibinfo{volume}{19}
	(\bibinfo{year}{2016}) \bibinfo{pages}{1172--1185}.
	\bibitem[{Bullmore and Sporns(2012)}]{bullmore2012economy}
	\bibinfo{author}{E.~Bullmore}, \bibinfo{author}{O.~Sporns},
	\newblock \bibinfo{title}{The economy of brain network organization},
	\newblock \bibinfo{journal}{Nature Reviews Neuroscience} \bibinfo{volume}{13}
	(\bibinfo{year}{2012}) \bibinfo{pages}{336--349}.
	\bibitem[{Ji et~al.(2023)Ji, Ye, Mu, Lin, Tian, Hens, Perc, Tang, Sun, and
		Kurths}]{ji2023signal}
	\bibinfo{author}{P.~Ji}, \bibinfo{author}{J.~Ye}, \bibinfo{author}{Y.~Mu},
	\bibinfo{author}{W.~Lin}, \bibinfo{author}{Y.~Tian},
	\bibinfo{author}{C.~Hens}, \bibinfo{author}{M.~Perc},
	\bibinfo{author}{Y.~Tang}, \bibinfo{author}{J.~Sun},
	\bibinfo{author}{J.~Kurths},
	\newblock \bibinfo{title}{Signal propagation in complex networks},
	\newblock \bibinfo{journal}{Physics Reports} \bibinfo{volume}{1017}
	(\bibinfo{year}{2023}) \bibinfo{pages}{1--96}.
	\bibitem[{Barabasi and Oltvai(2004)}]{barabasi2004network}
	\bibinfo{author}{A.-L. Barabasi}, \bibinfo{author}{Z.~N. Oltvai},
	\newblock \bibinfo{title}{Network biology: understanding the cell's functional
		organization},
	\newblock \bibinfo{journal}{Nature Reviews Genetics} \bibinfo{volume}{5}
	(\bibinfo{year}{2004}) \bibinfo{pages}{101--113}.
	\bibitem[{Shekatkar et~al.(2015)Shekatkar, Bhagwat, and
		Ambika}]{shekatkar2015divisibility}
	\bibinfo{author}{S.~M. Shekatkar}, \bibinfo{author}{C.~Bhagwat},
	\bibinfo{author}{G.~Ambika},
	\newblock \bibinfo{title}{Divisibility patterns of natural numbers on a complex
		network},
	\newblock \bibinfo{journal}{Scientific Reports} \bibinfo{volume}{5}
	(\bibinfo{year}{2015}) \bibinfo{pages}{1--10}.
	\bibitem[{Garcia-Perez et~al.(2014)Garcia-Perez, Serrano, and
		Bogun{\'a}}]{garcia2014complex}
	\bibinfo{author}{G.~Garcia-Perez}, \bibinfo{author}{M.~{\'A}. Serrano},
	\bibinfo{author}{M.~Bogun{\'a}},
	\newblock \bibinfo{title}{Complex architecture of primes and natural numbers},
	\newblock \bibinfo{journal}{Physical Review E} \bibinfo{volume}{90}
	(\bibinfo{year}{2014}) \bibinfo{pages}{022806}.
	\bibitem[{Solares-Hern{\'a}ndez et~al.(2020)Solares-Hern{\'a}ndez, Manzano,
		P{\'e}rez-Benito, and Conejero}]{solares2020divisibility}
	\bibinfo{author}{P.~A. Solares-Hern{\'a}ndez}, \bibinfo{author}{F.~A. Manzano},
	\bibinfo{author}{F.~J. P{\'e}rez-Benito}, \bibinfo{author}{J.~A. Conejero},
	\newblock \bibinfo{title}{Divisibility patterns within pascal divisibility
		networks},
	\newblock \bibinfo{journal}{Mathematics} \bibinfo{volume}{8}
	(\bibinfo{year}{2020}) \bibinfo{pages}{254}.
	\bibitem[{Chandra and Dasgupta(2005)}]{chandra2005small}
	\bibinfo{author}{A.~K. Chandra}, \bibinfo{author}{S.~Dasgupta},
	\newblock \bibinfo{title}{A small world network of prime numbers},
	\newblock \bibinfo{journal}{Physica A: Statistical Mechanics and its
		Applications} \bibinfo{volume}{357} (\bibinfo{year}{2005})
	\bibinfo{pages}{436--446}.
	\bibitem[{Thomason(1987)}]{thomason1987pseudo}
	\bibinfo{author}{A.~Thomason},
	\newblock \bibinfo{title}{Pseudo-random graphs},
	\newblock in: \bibinfo{booktitle}{North-Holland Mathematics Studies}, volume
	\bibinfo{volume}{144}, \bibinfo{publisher}{Elsevier}, \bibinfo{year}{1987},
	pp. \bibinfo{pages}{307--331}.
	\bibitem[{Krivelevich and Sudakov(2006)}]{krivelevich2006pseudo}
	\bibinfo{author}{M.~Krivelevich}, \bibinfo{author}{B.~Sudakov},
	\newblock \bibinfo{title}{Pseudo-random graphs},
	\newblock \bibinfo{journal}{More sets, graphs and numbers: A Salute to Vera Sos
		and Andr{\'a}s Hajnal}  (\bibinfo{year}{2006}) \bibinfo{pages}{199--262}.
	\bibitem[{Alon and Milman(1985)}]{alon1985lambda1}
	\bibinfo{author}{N.~Alon}, \bibinfo{author}{V.~D. Milman},
	\newblock \bibinfo{title}{$\lambda$1, isoperimetric inequalities for graphs,
		and superconcentrators},
	\newblock \bibinfo{journal}{Journal of Combinatorial Theory, Series B}
	\bibinfo{volume}{38} (\bibinfo{year}{1985}) \bibinfo{pages}{73--88}.
	\bibitem[{Alon(1986)}]{alon1986eigenvalues}
	\bibinfo{author}{N.~Alon},
	\newblock \bibinfo{title}{Eigenvalues and expanders},
	\newblock \bibinfo{journal}{Combinatorica} \bibinfo{volume}{6}
	(\bibinfo{year}{1986}) \bibinfo{pages}{83--96}.
	\bibitem[{Alon and Spencer(2016)}]{alon2016probabilistic}
	\bibinfo{author}{N.~Alon}, \bibinfo{author}{J.~H. Spencer}, \bibinfo{title}{The
		probabilistic method}, \bibinfo{publisher}{John Wiley \& Sons},
	\bibinfo{year}{2016}.
	\bibitem[{K{\"u}hn and Osthus(2006)}]{kuhn2006multicolored}
	\bibinfo{author}{D.~K{\"u}hn}, \bibinfo{author}{D.~Osthus},
	\newblock \bibinfo{title}{Multicolored hamilton cycles and perfect matchings in
		pseudorandom graphs},
	\newblock \bibinfo{journal}{SIAM Journal on Discrete Mathematics}
	\bibinfo{volume}{20} (\bibinfo{year}{2006}) \bibinfo{pages}{273--286}.
	\bibitem[{Krivelevich et~al.(2010)Krivelevich, Lee, and
		Sudakov}]{krivelevich2010resilient}
	\bibinfo{author}{M.~Krivelevich}, \bibinfo{author}{C.~Lee},
	\bibinfo{author}{B.~Sudakov},
	\newblock \bibinfo{title}{Resilient pancyclicity of random and pseudorandom
		graphs},
	\newblock \bibinfo{journal}{SIAM Journal on Discrete Mathematics}
	\bibinfo{volume}{24} (\bibinfo{year}{2010}) \bibinfo{pages}{1--16}.
	\bibitem[{Chung et~al.(1989)Chung, Graham, and Wilson}]{chung1989quasi}
	\bibinfo{author}{F.~R.~K. Chung}, \bibinfo{author}{R.~L. Graham},
	\bibinfo{author}{R.~M. Wilson},
	\newblock \bibinfo{title}{Quasi-random graphs},
	\newblock \bibinfo{journal}{Combinatorica} \bibinfo{volume}{9}
	(\bibinfo{year}{1989}) \bibinfo{pages}{345--362}.
	\bibitem[{Barab{\'a}si(2013)}]{barabasi2013network}
	\bibinfo{author}{A.-L. Barab{\'a}si},
	\newblock \bibinfo{title}{Network science},
	\newblock \bibinfo{journal}{Philosophical Transactions of the Royal Society A:
		Mathematical, Physical and Engineering Sciences} \bibinfo{volume}{371}
	(\bibinfo{year}{2013}) \bibinfo{pages}{20120375}.
	\bibitem[{Barahona and Pecora(2002)}]{barahona2002synchronization}
	\bibinfo{author}{M.~Barahona}, \bibinfo{author}{L.~M. Pecora},
	\newblock \bibinfo{title}{Synchronization in small-world systems},
	\newblock \bibinfo{journal}{Physical Review Letters} \bibinfo{volume}{89}
	(\bibinfo{year}{2002}) \bibinfo{pages}{054101}.
	\bibitem[{Pecora and Carroll(1998)}]{pecora1998master}
	\bibinfo{author}{L.~M. Pecora}, \bibinfo{author}{T.~L. Carroll},
	\newblock \bibinfo{title}{Master stability functions for synchronized coupled
		systems},
	\newblock \bibinfo{journal}{Physical Review Letters} \bibinfo{volume}{80}
	(\bibinfo{year}{1998}) \bibinfo{pages}{2109}.
	\bibitem[{Nishikawa et~al.(2003)Nishikawa, Motter, Lai, and
		Hoppensteadt}]{nishikawa2003heterogeneity}
	\bibinfo{author}{T.~Nishikawa}, \bibinfo{author}{A.~E. Motter},
	\bibinfo{author}{Y.-C. Lai}, \bibinfo{author}{F.~C. Hoppensteadt},
	\newblock \bibinfo{title}{Heterogeneity in oscillator networks: Are smaller
		worlds easier to synchronize?},
	\newblock \bibinfo{journal}{Physical Review Letters} \bibinfo{volume}{91}
	(\bibinfo{year}{2003}) \bibinfo{pages}{014101}.
	\bibitem[{Osipov et~al.(2007)Osipov, Kurths, and
		Zhou}]{osipov2007synchronization}
	\bibinfo{author}{G.~V. Osipov}, \bibinfo{author}{J.~Kurths},
	\bibinfo{author}{C.~Zhou}, \bibinfo{title}{Synchronization in oscillatory
		networks}, \bibinfo{publisher}{Springer Science \& Business Media},
	\bibinfo{year}{2007}.
	\bibitem[{Pikovsky et~al.(2001)Pikovsky, Rosenblum, and
		Kurths}]{pikovsky2001synchronization}
	\bibinfo{author}{A.~Pikovsky}, \bibinfo{author}{M.~Rosenblum},
	\bibinfo{author}{J.~Kurths},
	\newblock \bibinfo{title}{Synchronization cambridge university press},
	\newblock \bibinfo{journal}{Cambridge, England}  (\bibinfo{year}{2001}).
	\bibitem[{Dai et~al.(2020)Dai, Li, Guo, Jia, Perc, Manshour, Wang, and
		Boccaletti}]{dai2020discontinuous}
	\bibinfo{author}{X.~Dai}, \bibinfo{author}{X.~Li}, \bibinfo{author}{H.~Guo},
	\bibinfo{author}{D.~Jia}, \bibinfo{author}{M.~Perc},
	\bibinfo{author}{P.~Manshour}, \bibinfo{author}{Z.~Wang},
	\bibinfo{author}{S.~Boccaletti},
	\newblock \bibinfo{title}{Discontinuous transitions and rhythmic states in the
		d-dimensional kuramoto model induced by a positive feedback with the global
		order parameter},
	\newblock \bibinfo{journal}{Physical Review Letters} \bibinfo{volume}{125}
	(\bibinfo{year}{2020}) \bibinfo{pages}{194101}.
	\bibitem[{Apostol(1998)}]{apostol1998introduction}
	\bibinfo{author}{T.~M. Apostol}, \bibinfo{title}{Introduction to analytic
		number theory}, \bibinfo{publisher}{Springer Science \& Business Media},
	\bibinfo{year}{1998}.
	\bibitem[{T{\'o}th(2002)}]{toth2002probability}
	\bibinfo{author}{L.~T{\'o}th},
	\newblock \bibinfo{title}{The probability that k positive integers are pairwise
		relatively prime},
	\newblock \bibinfo{journal}{Fibonacci Quart} \bibinfo{volume}{40}
	(\bibinfo{year}{2002}) \bibinfo{pages}{13--18}.
	\bibitem[{Jafarizadeh et~al.(2020)Jafarizadeh, Tofigh, Lipman, and
		Abolhasan}]{jafarizadeh2020optimizing}
	\bibinfo{author}{S.~Jafarizadeh}, \bibinfo{author}{F.~Tofigh},
	\bibinfo{author}{J.~Lipman}, \bibinfo{author}{M.~Abolhasan},
	\newblock \bibinfo{title}{Optimizing synchronizability in networks of coupled
		systems},
	\newblock \bibinfo{journal}{Automatica} \bibinfo{volume}{112}
	(\bibinfo{year}{2020}) \bibinfo{pages}{108711}.
	\bibitem[{Tang et~al.(2014)Tang, Qian, Gao, and
		Kurths}]{tang2014synchronization}
	\bibinfo{author}{Y.~Tang}, \bibinfo{author}{F.~Qian}, \bibinfo{author}{H.~Gao},
	\bibinfo{author}{J.~Kurths},
	\newblock \bibinfo{title}{Synchronization in complex networks and its
		application--a survey of recent advances and challenges},
	\newblock \bibinfo{journal}{Annual Reviews in Control} \bibinfo{volume}{38}
	(\bibinfo{year}{2014}) \bibinfo{pages}{184--198}.
	\bibitem[{Holland and Hastings(2008)}]{holland2008strong}
	\bibinfo{author}{M.~D. Holland}, \bibinfo{author}{A.~Hastings},
	\newblock \bibinfo{title}{Strong effect of dispersal network structure on
		ecological dynamics},
	\newblock \bibinfo{journal}{Nature} \bibinfo{volume}{456}
	(\bibinfo{year}{2008}) \bibinfo{pages}{792--794}.
	\bibitem[{Hardy et~al.(1979)Hardy, Wright et~al.}]{hardy1979introduction}
	\bibinfo{author}{G.~H. Hardy}, \bibinfo{author}{E.~M. Wright}, et~al.,
	\bibinfo{title}{An introduction to the theory of numbers},
	\bibinfo{publisher}{Oxford university press}, \bibinfo{year}{1979}.
	\bibitem[{Sondow(2009)}]{sondow2009ramanujan}
	\bibinfo{author}{J.~Sondow},
	\newblock \bibinfo{title}{Ramanujan primes and bertrand's postulate},
	\newblock \bibinfo{journal}{The American Mathematical Monthly}
	\bibinfo{volume}{116} (\bibinfo{year}{2009}) \bibinfo{pages}{630--635}.
	
\end{thebibliography}

\end{document}